\documentclass[12pt]{article}%
\usepackage{amsfonts}
\usepackage{sw20bams}
\usepackage{amsmath}
\usepackage{amssymb}
\usepackage{graphicx}%
\providecommand{\U}[1]{\protect\rule{.1in}{.1in}}
\begin{document}

\title{Exercises in Iterational Asymptotics III}
\author{Steven Finch}
\date{March 17, 2025}
\maketitle

\begin{abstract}
The nonlinear recurrences we consider here include simple continued fractions
for the Golden \&\ Silver means and a parametric family of cubics in
connection with Abel's functional equation.

\end{abstract}

\footnotetext{Copyright \copyright \ 2025 by Steven R. Finch. All rights
reserved.}

This paper is a continuation of \cite{F1-exc3, F2-exc3}. \ The rate at which
simple continued fractions converge to familiar constants $\varphi,\psi$ turns
out to explicitly involve $\varphi,\psi$. \ Also, graphs of principal
solutions $F_{a}$ of $F_{a}(x-a\,x^{2}+x^{3})=F_{a}(x)+1$ for increasing
$a\in\lbrack1,2)$ exhibit behavior loosely resembling the topologist's sine curve.

\section{Huiti\`{e}me exercice}

\textbf{Consider initially} the famous recurrence
\[%
\begin{array}
[c]{ccccc}%
x_{k}=1+\dfrac{1}{x_{k-1}} &  & \text{for }k\geq1\text{;} &  & x_{0}=1.
\end{array}
\]
Quantify the convergence rate of $x_{k}$ as $k\rightarrow\infty$.

It is well known that%
\[
x_{0}=1<\dfrac{3}{2}=x_{2}<\ldots<\varphi<\ldots<x_{3}=\dfrac{5}{3}<2=x_{1}%
\]
where the limiting value%
\[
\varphi=\dfrac{1+\sqrt{5}}{2}=1.6180339887498948482045868...
\]
is the Golden mean. \ Letting
\[%
\begin{array}
[c]{ccc}%
f(x)=\dfrac{\varphi+x}{\varphi+1+x}-(\varphi-1), &  & g(x)=(\varphi
-1)-\dfrac{\varphi-x}{\varphi+1-x}%
\end{array}
\]
we have $x_{3}-\varphi=f(x_{1}-\varphi)$ and $\varphi-x_{2}=g(\varphi-x_{0})$.
\ The pattern is clear. \ Define%
\[%
\begin{array}
[c]{ccc}%
u_{k}=x_{2k+1}-\varphi, &  & v_{k}=\varphi-x_{2k}%
\end{array}
\]
and thus%
\[%
\begin{array}
[c]{ccc}%
u_{k+1}=f(u_{k}), &  & v_{k+1}=g(v_{k})
\end{array}
\]
for all $k$. \ Both $u_{k}$ and $v_{k}$ approach $0$; we determine the
respective speeds at which they do so, following Theorem 2.1 in
\cite{Thrn-exc3}. \ Note that $f(0)=g(0)=0$, $0<\max\{f(x),g(x)\}<x$ for all
$x>0$, and%
\[
f^{\prime}(0)=g^{\prime}(0)=\dfrac{1}{(1+\varphi)^{2}}=5-3\varphi
=0.1458980337503154553862395...<1.
\]
We now treat $f(x)$ and $g(x)$ separately.

The function%

\[
F(x)=\left\{
\begin{array}
[c]{lll}%
\dfrac{f(x)-\left(  5-3\varphi\right)  x}{x^{2}} &  & \text{if }x>0,\medskip\\
-\left(  13-8\varphi\right)  &  & \text{if }x=0
\end{array}
\right.
\]
is continuous and bounded on $[0,\infty)$; in fact, $\left\vert
F(x)\right\vert \leq M=13-8\varphi$ (left endpoint). \ Since $u_{k}$ is
decreasing, $u_{k}\leq u_{0}=2-\varphi$ always. \ Because%
\begin{align*}
\frac{u_{k+1}}{u_{k}}  &  =\frac{f(u_{k})}{u_{k}}=\left(  5-3\varphi\right)
+F(u_{k})\,u_{k}\leq\left(  5-3\varphi\right)  +\left(  13-8\varphi\right)
\left(  2-\varphi\right)  =39-24\varphi\\
&  <39-24\cdot\frac{8}{5}=\frac{3}{5}%
\end{align*}
we have%
\[
u_{k+1}<\left(  \frac{3}{5}\right)  u_{k}<\left(  \frac{3}{5}\right)
^{2}u_{k-1}<\ldots<\left(  \frac{3}{5}\right)  ^{k+1}u_{0}.
\]
It follows that the series%
\[
(1+\varphi)^{2}\,%
%TCIMACRO{\dsum \limits_{k=0}^{\infty}}%
%BeginExpansion
{\displaystyle\sum\limits_{k=0}^{\infty}}
%EndExpansion
\,u_{k}\left\vert F(u_{k})\right\vert <(1+\varphi)^{2}M\,%
%TCIMACRO{\dsum \limits_{k=0}^{\infty}}%
%BeginExpansion
{\displaystyle\sum\limits_{k=0}^{\infty}}
%EndExpansion
\,u_{k}<(1+\varphi)^{2}M\,u_{0}\,%
%TCIMACRO{\dsum \limits_{k=0}^{\infty}}%
%BeginExpansion
{\displaystyle\sum\limits_{k=0}^{\infty}}
%EndExpansion
\left(  \frac{3}{5}\right)  ^{k}%
\]
converges, which in turn implies that the product%
\[%
%TCIMACRO{\dprod \limits_{k=0}^{\infty}}%
%BeginExpansion
{\displaystyle\prod\limits_{k=0}^{\infty}}
%EndExpansion
\left(  1+(1+\varphi)^{2}\,u_{k}\,F(u_{k})\right)
\]
also converges. \ Finally, multiplying both sides of
\[
(1+\varphi)^{2}\frac{u_{j+1}}{u_{j}}=1+(1+\varphi)^{2}\,u_{j}\,F(u_{j})
\]
from $j=0$ to $k-1$ gives%
\[
(1+\varphi)^{2k}\,\frac{u_{k}}{u_{0}}=%
%TCIMACRO{\dprod \limits_{j=0}^{k-1}}%
%BeginExpansion
{\displaystyle\prod\limits_{j=0}^{k-1}}
%EndExpansion
\left(  1+(1+\varphi)^{2}\,u_{j}\,F(u_{j})\right)
\]
and therefore%
\begin{align*}
\lim_{k\rightarrow\infty}(1+\varphi)^{2k}\,u_{k}  &  =(2-\varphi)%
%TCIMACRO{\dprod \limits_{j=0}^{\infty}}%
%BeginExpansion
{\displaystyle\prod\limits_{j=0}^{\infty}}
%EndExpansion
\left(  1+(1+\varphi)^{2}\,u_{j}\,F(u_{j})\right)
=0.3262379212492639374321078...\\
&  =7\varphi-11.
\end{align*}
Having finished with $f(x)$, we now investigate $g(x)$.

The function%

\[
G(x)=\left\{
\begin{array}
[c]{lll}%
\dfrac{g(x)-\left(  5-3\varphi\right)  x}{x^{2}} &  & \text{if }x>0,\medskip\\
13-8\varphi &  & \text{if }x=0
\end{array}
\right.
\]
is continuous and bounded on $[0,\varphi-1]$; in fact, $\left\vert
G(x)\right\vert \leq M=\left(  5-3\varphi\right)  /2$ (right endpoint).
\ \ Since $v_{k}$ is decreasing, $v_{k}\leq v_{0}=\varphi-1$ always. \ A
similar line of reasoning gives%
\begin{align*}
\lim_{k\rightarrow\infty}(1+\varphi)^{2k}\,v_{k}  &  =(\varphi-1)%
%TCIMACRO{\dprod \limits_{j=0}^{\infty}}%
%BeginExpansion
{\displaystyle\prod\limits_{j=0}^{\infty}}
%EndExpansion
\left(  1+(1+\varphi)^{2}\,v_{j}\,G(v_{j})\right)
=0.8541019662496845446137605...\\
&  =3\varphi-4.
\end{align*}
The two constants here differ by a factor of $1+\varphi$. \ 

\textbf{Consider finally} the recurrence
\[%
\begin{array}
[c]{ccccc}%
x_{k}=2+\dfrac{1}{x_{k-1}} &  & \text{for }k\geq1\text{;} &  & x_{0}=2.
\end{array}
\]
Again, quantify the convergence rate of $x_{k}$ as $k\rightarrow\infty$.

Here we have%
\[
x_{0}=2<\dfrac{12}{5}=x_{2}<\ldots<\psi<\ldots<x_{3}=\dfrac{29}{12}<\dfrac
{5}{2}=x_{1}%
\]
where the limiting value%
\[
\psi=1+\sqrt{2}=2.4142135623730950488016887...
\]
is the Silver mean. \ Letting
\[%
\begin{array}
[c]{ccc}%
f(x)=\dfrac{\psi+x}{2\psi+1+2x}-(\psi-2), &  & g(x)=(\psi-2)-\dfrac{\psi
-x}{2\psi+1-2x}%
\end{array}
\]
we have $x_{3}-\psi=f(x_{1}-\psi)$, $\psi-x_{2}=g(\psi-x_{0})$ and%
\[
f^{\prime}(0)=g^{\prime}(0)=\dfrac{1}{(1+2\psi)^{2}}=29-12\psi
=0.0294372515228594143797353...<1.
\]
Defining $u_{k}$ and $v_{k}$ analogously and omitting details, it follows that%
\[
F(x)=\left\{
\begin{array}
[c]{lll}%
\dfrac{f(x)-\left(  29-12\psi\right)  x}{x^{2}} &  & \text{if }x>0,\medskip\\
-\left(  338-140\psi\right)  &  & \text{if }x=0
\end{array}
\right.
\]%
\begin{align*}
\lim_{k\rightarrow\infty}(1+2\psi)^{2k}\,u_{k}  &  =\left(  \dfrac{5}{2}%
-\psi\right)
%TCIMACRO{\dprod \limits_{j=0}^{\infty}}%
%BeginExpansion
{\displaystyle\prod\limits_{j=0}^{\infty}}
%EndExpansion
\left(  1+(1+2\psi)^{2}\,u_{j}\,F(u_{j})\right)
=0.0832611206852316592574166...\\
&  =34\psi-82
\end{align*}
and%
\[
G(x)=\left\{
\begin{array}
[c]{lll}%
\dfrac{g(x)-\left(  29-12\psi\right)  x}{x^{2}} &  & \text{if }x>0,\medskip\\
338-140\psi &  & \text{if }x=0
\end{array}
\right.
\]%
\begin{align*}
\lim_{k\rightarrow\infty}(1+2\psi)^{2k}\,v_{k}  &  =(\psi-2)%
%TCIMACRO{\dprod \limits_{j=0}^{\infty}}%
%BeginExpansion
{\displaystyle\prod\limits_{j=0}^{\infty}}
%EndExpansion
\left(  1+(1+2\psi)^{2}\,v_{j}\,G(v_{j})\right)
=0.4852813742385702928101323...\\
&  =6\psi-14.
\end{align*}
The two constants here differ by a factor of $1+2\psi$. \ 

From a different point of view (rational approximation, Hurwitz's theorem and
the Lagrange spectrum), the Golden mean is \textquotedblleft
worse\textquotedblright\ than the Silver mean \cite{F0-exc3}. \ This is
consistent with our observation that $1+\varphi\approx2.618<5.828\approx
1+2\psi$, i.e., convergence (in our sense) is slower for continued fractions
with all $1$s (as partial denominators) than those with all $2$s.

\section{Neuvi\`{e}me exercice}

\textbf{Let the function} $f_{a}:[0,a]\rightarrow\lbrack0,a]$ be defined by
$f_{a}(x)=x\left(  1-a\,x+x^{2}\right)  $, where $1\leq a\leq2$ is a
parameter. \ Let $F_{a}$ denote the solution, via $f_{a}$ iterates, of Abel's
equation $F_{a}(f_{a}(x))=F_{a}(x)+1$. \ Describe what happens graphically as
$a$ increases. \ %

%TCIMACRO{\FRAME{ftbpFU}{2.8003in}{2.8496in}{0pt}{\Qcb{Blue curve is
%$10f_{a}(x)$, scaled for visibility; green curve is $F_{a}(x)$; parameter
%$a=3/2$. \ $f_{a}$ is increasing; $F_{a}$ is decreasing. \ Distance between
%vertical notches is four times that for horizontal.}}{}{ex1.eps}%
%width 2.8003in;  height 2.8496in;  depth 0pt;  original-width 8.0125in;
%original-height 8.1552in;  cropleft "0";  croptop "1";  cropright "1";
%cropbottom "0";  filename 'ex1.eps';file-properties "XNPEU";}} }%
%BeginExpansion
\begin{figure}
[ptb]
\begin{center}
\includegraphics[
height=2.8496in,
width=2.8003in
]%
{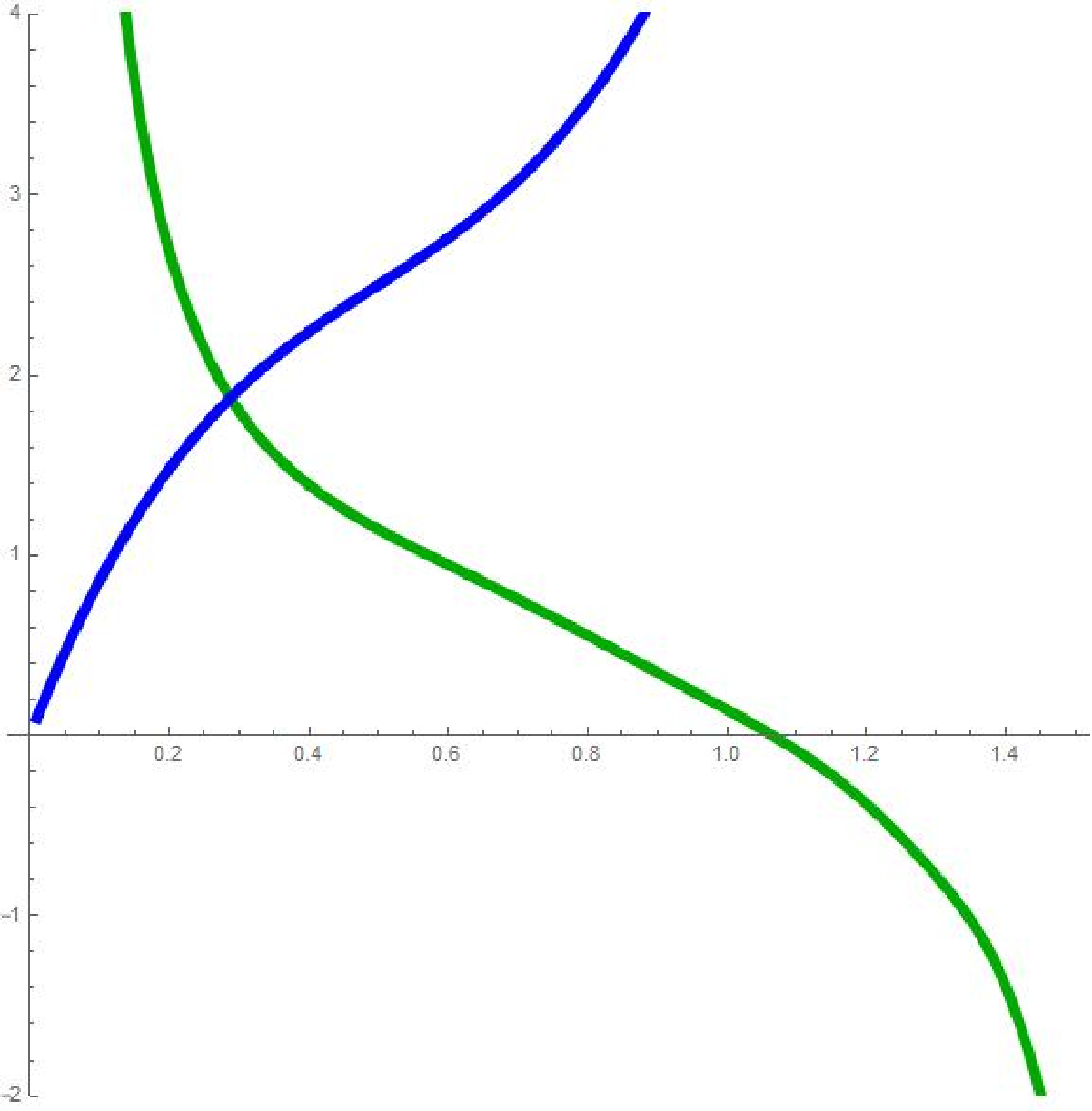}%
\caption{Blue curve is $10f_{a}(x)$, scaled for visibility; green curve is
$F_{a}(x)$; parameter $a=3/2$. \ $f_{a}$ is increasing; $F_{a}$ is decreasing.
\ Distance between vertical notches is four times that for horizontal.}%
\end{center}
\end{figure}
%EndExpansion
%TCIMACRO{\FRAME{ftbpFU}{2.8011in}{2.8504in}{0pt}{\Qcb{Parameter $a=5/3$.
%\ Both curves remain monotone; gentle undulations start to appear in $F_{a}$.
%\ }}{}{ex2.eps}{\special{ language "Scientific Word";  type "GRAPHIC";
%width 2.8011in;  height 2.8504in;  depth 0pt;  original-width 8.0955in;
%original-height 8.2382in;  cropleft "0";  croptop "1";  cropright "1";
%cropbottom "0";  filename 'ex2.eps';file-properties "XNPEU";}} }%
%BeginExpansion
\begin{figure}
[ptb]
\begin{center}
\includegraphics[
height=2.8504in,
width=2.8011in
]%
{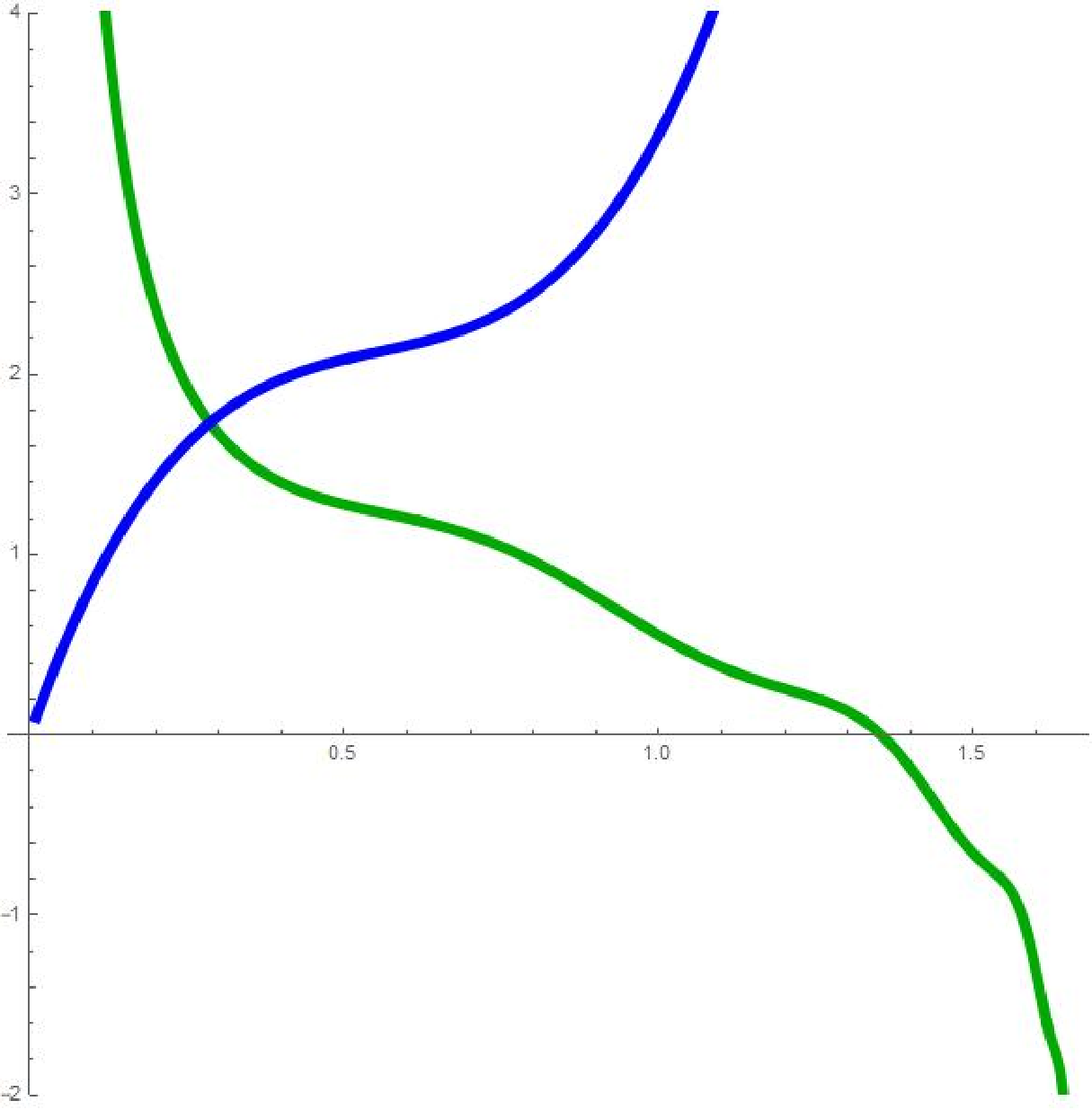}%
\caption{Parameter $a=5/3$. \ Both curves remain monotone; gentle undulations
start to appear in $F_{a}$. \ }%
\end{center}
\end{figure}
%EndExpansion
%TCIMACRO{\FRAME{ftbpFU}{2.8124in}{2.8496in}{0pt}{\Qcb{Parameter $a=\sqrt{3}$.
%Undulations become more pronounced. \ $F_{a}$-critical points are at
%$x=0.577,1.304,1.613,1.701,\ldots.$}}{}{ex3.eps}%
%width 2.8124in;  height 2.8496in;  depth 0pt;  original-width 8.4034in;
%original-height 8.5167in;  cropleft "0";  croptop "1";  cropright "1";
%cropbottom "0";  filename 'ex3.eps';file-properties "XNPEU";}} }%
%BeginExpansion
\begin{figure}
[ptb]
\begin{center}
\includegraphics[
height=2.8496in,
width=2.8124in
]%
{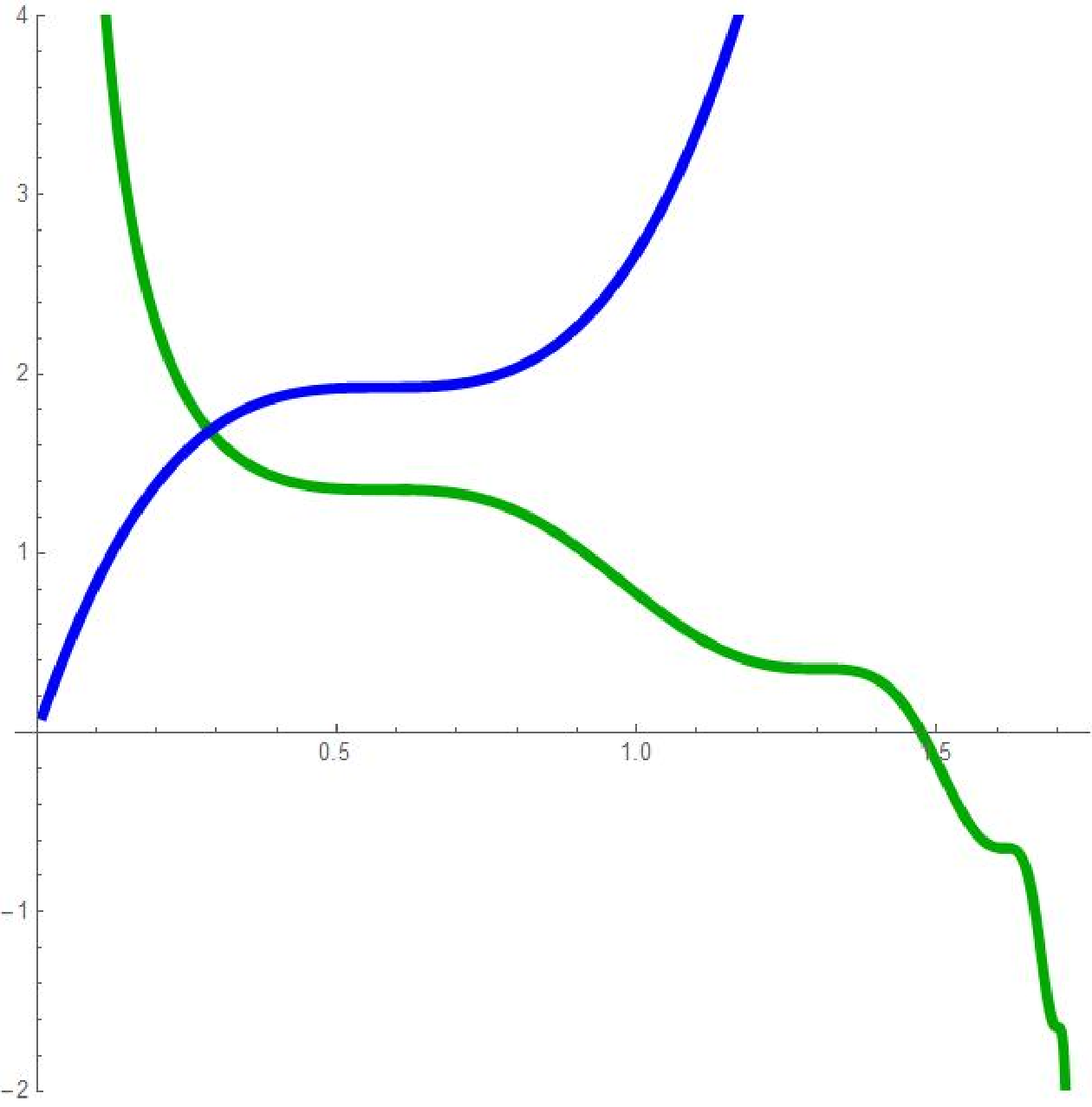}%
\caption{Parameter $a=\sqrt{3}$. Undulations become more pronounced. \ $F_{a}%
$-critical points are at $x=0.577,1.304,1.613,1.701,\ldots.$}%
\end{center}
\end{figure}
%EndExpansion
%TCIMACRO{\FRAME{ftbpFU}{2.8029in}{2.8496in}{0pt}{\Qcb{Parameter $a=9/5$.
%Curves are no longer monotone. \ $F_{a}$-minimum points appear at
%$x=0.436,1.286,1.663,1.766\ldots$; $F_{a}$-maximum points appear at
%$x=0.763,1.472,1.717,1.780\ldots.$}}{}{ex4.eps}%
%width 2.8029in;  height 2.8496in;  depth 0pt;  original-width 8.3454in;
%original-height 8.489in;  cropleft "0";  croptop "1";  cropright "1";
%cropbottom "0";  filename 'ex4.eps';file-properties "XNPEU";}} }%
%BeginExpansion
\begin{figure}
[ptb]
\begin{center}
\includegraphics[
height=2.8496in,
width=2.8029in
]%
{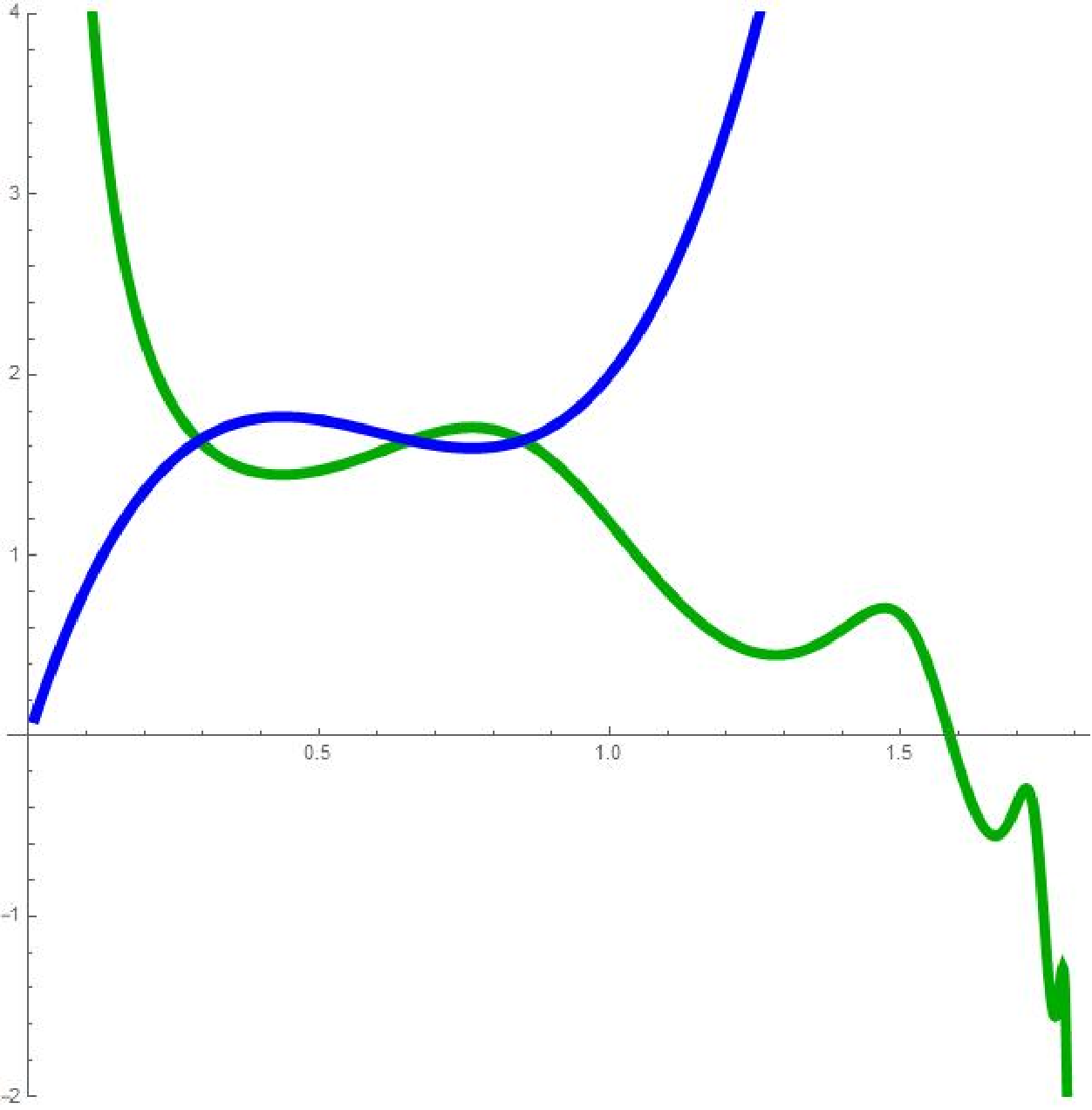}%
\caption{Parameter $a=9/5$. Curves are no longer monotone. \ $F_{a}$-minimum
points appear at $x=0.436,1.286,1.663,1.766\ldots$; $F_{a}$-maximum points
appear at $x=0.763,1.472,1.717,1.780\ldots.$}%
\end{center}
\end{figure}
%EndExpansion
%TCIMACRO{\FRAME{ftbpFU}{2.8029in}{2.8496in}{0pt}{\Qcb{Parameter $a=19/10$.
%Violent waves dominate the graph. \ $F_{a}$-minimum points appear at
%$x=0.372,1.368,1.771,1.871\ldots$; $F_{a}$-maximum points appear at
%$x=0.893,1.623,1.836,1.886\ldots.$}}{}{ex5.eps}%
%width 2.8029in;  height 2.8496in;  depth 0pt;  original-width 8.4034in;
%original-height 8.5443in;  cropleft "0";  croptop "1";  cropright "1";
%cropbottom "0";  filename 'ex5.eps';file-properties "XNPEU";}} }%
%BeginExpansion
\begin{figure}
[ptb]
\begin{center}
\includegraphics[
height=2.8496in,
width=2.8029in
]%
{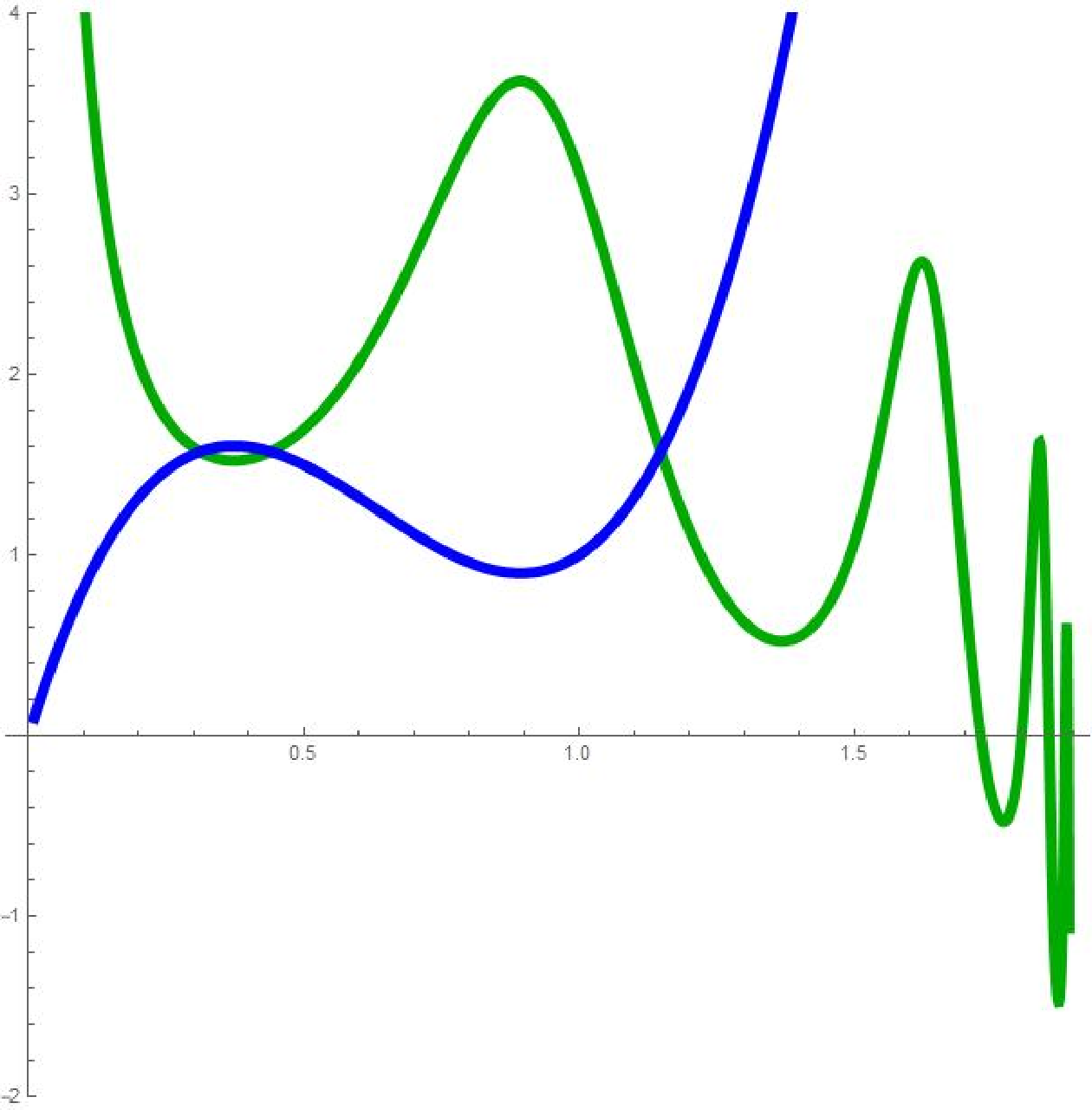}%
\caption{Parameter $a=19/10$. Violent waves dominate the graph. \ $F_{a}%
$-minimum points appear at $x=0.372,1.368,1.771,1.871\ldots$; $F_{a}$-maximum
points appear at $x=0.893,1.623,1.836,1.886\ldots.$}%
\end{center}
\end{figure}
%EndExpansion
%TCIMACRO{\FRAME{ftbpFU}{2.8029in}{2.8496in}{0pt}{\Qcb{Parameter $a=2$.
%\ Limiting scenario of the preceding. $F_{a}$-minimum points appear at
%$x=0.333,1.475,1.884\ldots$; $F_{a}$-vertical asymptotes appear at
%$x=1.000,1.754,1.948\ldots.$}}{}{ex6.eps}%
%width 2.8029in;  height 2.8496in;  depth 0pt;  original-width 8.4034in;
%original-height 8.5443in;  cropleft "0";  croptop "1";  cropright "1";
%cropbottom "0";  filename 'ex6.eps';file-properties "XNPEU";}} }%
%BeginExpansion
\begin{figure}
[ptb]
\begin{center}
\includegraphics[
height=2.8496in,
width=2.8029in
]%
{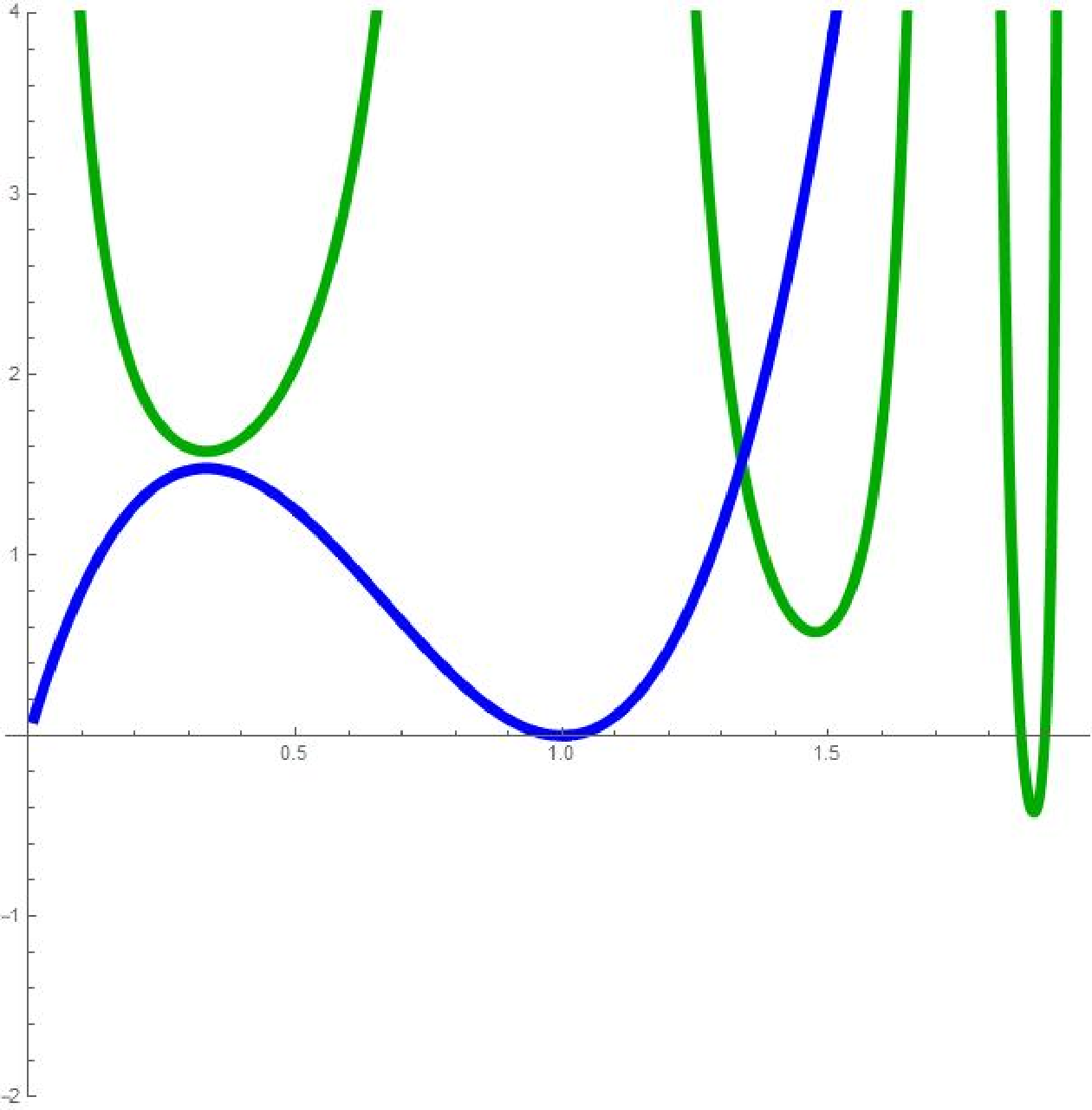}%
\caption{Parameter $a=2$. \ Limiting scenario of the preceding. $F_{a}%
$-minimum points appear at $x=0.333,1.475,1.884\ldots$; $F_{a}$-vertical
asymptotes appear at $x=1.000,1.754,1.948\ldots.$}%
\end{center}
\end{figure}
%EndExpansion

Let us focus on locations where derivatives vanish. \ Differentiating both
sides of the equation, we obtain%
\[
F_{a}^{\prime}(x-a\,x^{2}+x^{3})\left(  1-2a\,x+3x^{2}\right)  =F_{a}^{\prime
}(x)
\]
thus the leftmost pair of critical points comes from the quadratic factor:%
\[%
\begin{array}
[c]{ccc}%
\xi_{0}=\dfrac{a-\sqrt{a^{2}-3}}{3}, &  & \eta_{0}=\dfrac{a+\sqrt{a^{2}-3}}{3}%
\end{array}
\]
assuming $a\geq\sqrt{3}$. \ If $a>\sqrt{3}$, then $\xi_{0}$\ corresponds to
the leftmost minimum and $\eta_{0}$ corresponds to the leftmost maximum.
\ Moving forward, define $\xi_{k},\eta_{k}$ to be the $k^{\text{th}}$ pair of
critical points, i.e., unique real solutions of the cubic system%
\[%
\begin{array}
[c]{ccccc}%
x-a\,x^{2}+x^{3}=\xi_{k-1}, &  & y-a\,y^{2}+y^{3}=\eta_{k-1}, &  & k\geq1.
\end{array}
\]
For example, if $a=\sqrt{3}$, then $\xi_{k}=\eta_{k}$ always and%
\[%
\begin{array}
[c]{ccc}%
\xi_{0}=\dfrac{1}{\sqrt{3}}=0.577350269189..., &  & F(\xi_{0}%
)=1.354567323982...;
\end{array}
\]%
\[%
\begin{array}
[c]{ccc}%
\xi_{1}=1.304766026504..., &  & F(\xi_{1})=0.354567323982...;
\end{array}
\]%
\[%
\begin{array}
[c]{ccc}%
\xi_{2}=1.613468669954..., &  & F(\xi_{2})=-0.645432676017...;
\end{array}
\]%
\[%
\begin{array}
[c]{ccc}%
\xi_{3}=1.701609819539..., &  & F(\xi_{3})=-1.645432676017....
\end{array}
\]
If $a=9/5$, then
\[%
\begin{array}
[c]{ccc}%
\xi_{0}=\dfrac{9-\sqrt{6}}{15}=0.436700683814..., &  & F(\xi_{0}%
)=1.446679716680...;
\end{array}
\]%
\[%
\begin{array}
[c]{ccc}%
\xi_{1}=1.286564033401..., &  & F(\xi_{1})=0.446679716680...;
\end{array}
\]%
\[%
\begin{array}
[c]{ccc}%
\xi_{2}=1.663738765627..., &  & F(\xi_{2})=-0.553320283319...;
\end{array}
\]%
\[%
\begin{array}
[c]{ccc}%
\xi_{3}=1.766943652708..., &  & F(\xi_{3})=-1.553320283319...
\end{array}
\]
and%
\[%
\begin{array}
[c]{ccc}%
\eta_{0}=\dfrac{9+\sqrt{6}}{15}=0.763299316185..., &  & F(\eta_{0}%
)=1.707702719524...;
\end{array}
\]%
\[%
\begin{array}
[c]{ccc}%
\eta_{1}=1.472909661275..., &  & F(\eta_{1})=0.707702719524...;
\end{array}
\]%
\[%
\begin{array}
[c]{ccc}%
\eta_{2}=1.717164058718..., &  & F(\eta_{2})=-0.292297280475...;
\end{array}
\]%
\[%
\begin{array}
[c]{ccc}%
\eta_{3}=1.780129844535..., &  & F(\eta_{3})=-1.292297280475....
\end{array}
\]
If $a=19/10$, then
\[%
\begin{array}
[c]{ccc}%
\xi_{0}=\dfrac{19-\sqrt{61}}{30}=0.372991677469..., &  & F(\xi_{0}%
)=1.524207513358...;
\end{array}
\]%
\[%
\begin{array}
[c]{ccc}%
\xi_{1}=1.368415902116..., &  & F(\xi_{1})=0.524207513358...;
\end{array}
\]%
\[%
\begin{array}
[c]{ccc}%
\xi_{2}=1.771547833683..., &  & F(\xi_{2})=-0.475792486641...;
\end{array}
\]%
\[%
\begin{array}
[c]{ccc}%
\xi_{3}=1.871470296788..., &  & F(\xi_{3})=-1.475792486641...
\end{array}
\]
and%
\[%
\begin{array}
[c]{ccc}%
\eta_{0}=\dfrac{19+\sqrt{61}}{30}=0.893674989196..., &  & F(\eta
_{0})=3.626360576962...;
\end{array}
\]%
\[%
\begin{array}
[c]{ccc}%
\eta_{1}=1.623120175823..., &  & F(\eta_{1})=2.626360576962...;
\end{array}
\]%
\[%
\begin{array}
[c]{ccc}%
\eta_{2}=1.836690492410..., &  & F(\eta_{2})=1.626360576962...;
\end{array}
\]%
\[%
\begin{array}
[c]{ccc}%
\eta_{3}=1.886108430692..., &  & F(\eta_{3})=0.626360576962....
\end{array}
\]
If $a=2$, then
\[%
\begin{array}
[c]{ccc}%
\xi_{0}=\dfrac{1}{3}=0.333333333333..., &  & F(\xi_{0})=1.574672245867...;
\end{array}
\]%
\[%
\begin{array}
[c]{ccc}%
\xi_{1}=1.475329585787..., &  & F(\xi_{1})=0.574672245867...;
\end{array}
\]%
\[%
\begin{array}
[c]{ccc}%
\xi_{2}=1.884745179770..., &  & F(\xi_{2})=-0.425327754132...;
\end{array}
\]
and
\[%
\begin{array}
[c]{ccccc}%
\eta_{0}=1, &  & \eta_{1}=1.754877666246..., &  & \eta_{2}=1.948914407000....
\end{array}
\]
The latter are vertical asymptotes. \ Our analogy between $F_{a}(x)$ as
$a\rightarrow2^{-}$ and $\sin(1/(2-x))$ should not be interpreted too
literally. \ The similarity is only with horizontal compression in a
neighborhood of $x=2^{-}$; vertically $F_{a}$ explodes (while descending, in
essence) whereas $\sin$ remains bounded (and stationary, of course).

Computations here were performed using the highly precise Mavecha-Laohakosol
algorithm \cite{ML-exc3, F3-exc3} based on work by de Bruijn \cite{dB-exc3}
and Bencherif \&\ Robin \cite{BR-exc3}. \ Similar work was done in
\cite{F4-exc3}. \ The constants $\frac{1}{2}F_{2}\left(  \frac{1}{2}\right)  $
and $\frac{1}{2}F_{2}\left(  \frac{1}{3}\right)  $ also appear in
\cite{F5-exc3}.

\section{Dixi\`{e}me exercice}

\textbf{Given the recurrence} \cite{F4-exc3}%
\[%
\begin{array}
[c]{ccc}%
x_{n+1}=x_{n}+1+\dfrac{1}{x_{n}^{2}} &  & \text{for }n\geq0
\end{array}
\]
find the values $C(2)$, $C^{\prime}(2)$, $C^{\prime\prime}(2)$ in two distinct
ways, where $C(x_{0})=\lim_{n\rightarrow\infty}\left(  x_{n}-n\right)  $.

We infer from a similar example in \cite{IS-exc3} that%
\[
C(x_{0})=x_{0}+%
%TCIMACRO{\dsum \limits_{n=0}^{\infty}}%
%BeginExpansion
{\displaystyle\sum\limits_{n=0}^{\infty}}
%EndExpansion
\,\frac{1}{x_{n}^{2}}=2.598786...,
\]%
\[
C^{\prime}(x_{0})=%
%TCIMACRO{\dprod \limits_{n=0}^{\infty}}%
%BeginExpansion
{\displaystyle\prod\limits_{n=0}^{\infty}}
%EndExpansion
\left(  1-\frac{2}{x_{n}^{3}}\right)  =0.6615613240486...
\]
and, supposing the truth of a uniformity condition \cite{LTB-exc3},%
\[
C^{\prime\prime}(x_{0})=C^{\prime}(x_{0})\,%
%TCIMACRO{\dsum \limits_{n=0}^{\infty}}%
%BeginExpansion
{\displaystyle\sum\limits_{n=0}^{\infty}}
%EndExpansion
\left(  1-\frac{2}{x_{n}^{3}}\right)  ^{-1}\frac{6}{x_{n}^{4}}=0.3909...(?)
\]
In fact, the $0.3909$ estimate appears to be wrong and the series might simply
not be uniformly convergent (as a function of $x_{0}$). \ A more
computationally efficient approach makes use of the asymptotic expansion%

\begin{align*}
x_{n}  &  \sim n+C-\frac{1}{n}-\left(  \frac{1}{2}-C\right)  \frac{1}{n^{2}%
}+\left(  -\frac{5}{6}+C-C^{2}\right)  \frac{1}{n^{3}}+\left(  -\frac{5}%
{4}+\frac{5}{2}C-\frac{3}{2}C^{2}+C^{3}\right)  \frac{1}{n^{4}}\\
&  +\left(  -\frac{31}{15}+5C-5C^{2}+2C^{3}-C^{4}\right)  \frac{1}{n^{5}%
}+\left(  -\frac{11}{3}+\frac{31}{3}C-\frac{25}{2}C^{2}+\frac{25}{3}%
C^{3}-\frac{5}{2}C^{4}+C^{5}\right)  \frac{1}{n^{6}}\\
&  +\left(  -\frac{473}{70}+22C-31C^{2}+25C^{3}-\frac{25}{4}C^{4}+3C^{5}%
-C^{6}\right)  \frac{1}{n^{7}}\\
&  +\left(  -\frac{511}{40}+\frac{473}{10}C-77C^{2}+\frac{217}{3}C^{3}%
-\frac{175}{4}C^{4}+\frac{35}{2}C^{5}-\frac{7}{2}C^{6}+C^{7}\right)  \frac
{1}{n^{8}}\\
&  +\left(  -\frac{46651}{1890}+\frac{511}{5}C-\frac{946}{5}C^{2}+\frac
{616}{3}C^{3}-\frac{434}{3}C^{4}+70C^{5}-\frac{70}{3}C^{6}+4C^{7}%
-C^{8}\right)  \frac{1}{n^{9}}\\
&  +\left(  -\frac{20401}{420}+\frac{46651}{210}C-\frac{4599}{10}C^{2}%
+\frac{2838}{5}C^{3}-462C^{4}\right. \\
&  \left.  +\;\frac{1302}{5}C^{5}-105C^{6}+30C^{7}-\frac{9}{2}C^{8}%
+C^{9}\right)  \frac{1}{n^{10}}+\left(  -\frac{16124719}{166320}+\frac
{20401}{42}C-\frac{46651}{42}C^{2}\right. \\
&  \left.  +\;1533C^{3}-1419C^{4}+924C^{5}-434C^{6}+150C^{7}-\frac{75}{2}%
C^{8}+5C^{9}-C^{10}\right)  \frac{1}{n^{11}}%
\end{align*}%
\begin{align*}
&  +\left(  -\frac{1183139}{6048}+\frac{16124719}{15120}C-\frac{224411}%
{84}C^{2}+\frac{513161}{126}C^{3}-\frac{16863}{4}C^{4}+\frac{15609}{5}%
C^{5}\right. \\
&  \left.  -\;1694C^{6}+682C^{7}-\frac{825}{4}C^{8}+\frac{275}{6}C^{9}%
-\frac{11}{2}C^{10}+C^{11}\right)  \frac{1}{n^{12}}\\
&  +\left(  -\frac{1076978467}{2702700}+\frac{1183139}{504}C-\frac
{16124719}{2520}C^{2}+\frac{224411}{21}C^{3}-\frac{513161}{42}C^{4}%
+\frac{50589}{5}C^{5}\right. \\
&  \left.  -\;\frac{31218}{5}C^{6}+2904C^{7}-1023C^{8}+275C^{9}-55C^{10}%
+6C^{11}-C^{12}\right)  \frac{1}{n^{13}}%
\end{align*}
as $n\rightarrow\infty$. \ This was obtained via the algorithm in
\cite{ML-exc3}. \ For $N=10^{7}$, we calculate $x_{N}$ exactly via recursion,
and then set the value $x_{N}$ equal to our series and numerically solve for
$C$. \ It follows that%
\[
C(2)=2.5987868558248713482599664951883194762422902129186367437296275388853210...
\]
and, setting $\varepsilon=10^{-20}$,%
\[
C^{\prime}(2)\approx\frac{f(2+\varepsilon)-f(2-\varepsilon)}{2\varepsilon
}=0.6615613240486860705677502635...,
\]%
\[
C^{\prime\prime}(2)\approx\frac{f(2+\varepsilon)-2f(2)+f(2-\varepsilon
)}{\varepsilon^{2}}=0.37462642198301734111....
\]
The latter does not match the earlier $0.3909$ estimate, contradicting our
uniformity supposition. \ An independent confirmation of the $0.3746$ estimate
remains open.

\section{Acknowledgements}

The creators of Mathematica earn my gratitude every day:\ this paper could not
have otherwise been written.


\begin{thebibliography}{99}                                                                                               %


\bibitem {F1-exc3}S. R. Finch, Exercises in iterational asymptotics, arXiv:2411.16062.

\bibitem {F2-exc3}S. R. Finch, Exercises in iterational asymptotics II, arXiv:2501.06065.

\bibitem {Thrn-exc3}W. J. Thron, Sequences generated by iteration,
\textit{Trans. Amer. Math. Soc}. 96 (1960) 38--53; MR0117462.

\bibitem {F0-exc3}S. R. Finch, Freiman's constant, \textit{Mathematical
Constants}, Cambridge Univ. Press, 2003, pp. 199--203; MR2003519.

\bibitem {ML-exc3}S. Mavecha and V. Laohakosol, Asymptotic expansions of
iterates of some classical functions, \textit{Appl. Math. E-Notes} 13 (2013)
77--91; MR3121616; http://www.emis.de/journals/AMEN/2013/2013.htm.

\bibitem {F3-exc3}S. R. Finch, What do $\sin(x)$ and $\operatorname{arcsinh}%
(x)$ have in common? arXiv:2411.01591.

\bibitem {dB-exc3}N. G. de Bruijn, \textit{Asymptotic Methods in Analysis},
North-Holland, 1958; MR0099564 / Dover, 1981; MR0671583.

\bibitem {BR-exc3}F. Bencherif and G. Robin, Sur l'it\'{e}r\'{e} de $\sin(x)$,
\textit{Publ. Inst. Math. (Beograd)} 56(70) (1994) 41--53; MR1349068; http://eudml.org/doc/256122.

\bibitem {F4-exc3}S. R. Finch, Abel's functional equation and interrelations, arXiv:2503.00579.

\bibitem {F5-exc3}S. R. Finch, Generalized logistic maps and convergence, arXiv:24409.15175.

\bibitem {IS-exc3}E. Ionascu and P. Stanica, Effective asymptotics for some
nonlinear recurrences and almost doubly-exponential sequences, \textit{Acta
Math. Univ. Comenian.} 73 (2004) 75--87; MR2076045.

\bibitem {LTB-exc3}C. H. C. Little, K. L. Teo and B. van Brunt, \textit{An
Introduction to Infinite Products}, Springer-Verlag, 2022, pp. 84--87; MR4393582.%

\begin{tabular}
[c]{lll}
& Steven Finch & \\
& MIT Sloan School of Management & \\
& Cambridge, MA, USA & \\
& \textit{steven\_finch\_math@outlook.com} &
\end{tabular}

\end{thebibliography}
\end{document}